\documentclass[12pt]{article}
\usepackage{latexsym}
\usepackage{amssymb}
\usepackage{graphicx}
\usepackage{cite}

\newtheorem{Theorem}{Theorem}[part]
\newtheorem{Definition}{Definition}[part]
\newtheorem{Proposition}{Proposition}[part]

\newtheorem{Lemma}{Lemma}[part]
\newtheorem{Corollary}{Corollary}[part]

\newtheorem{Example}{Example}[part]

\newtheorem{Question}{Question}[part]

\def \ep{\hbox{ }\hfill$\Box$}

\begin{document}
\title{Introduction to Tensor Variational Inequalities}

\author{Yong Wang\thanks{Department of Mathematics, School of Science, Tianjin University, Tianjin 300072, P.R. China.
Email: wang\_yong@tju.edu.cn. This author's work was supported by the National Natural Science Foundation of China (Grant No. 71572125).}
\and
Zheng-Hai Huang\thanks{Corresponding author. Department of Mathematics, School of Science, Tianjin University, Tianjin 300072, P.R. China.
Email: huangzhenghai@tju.edu.cn. Tel: +86-22-27403615, Fax: +86-22-27403615. This author's work was supported by the National Natural Science Foundation of China (Grant No. 11431002).}
\and
Liqun Qi\thanks{Department of Applied Mathematics, The Hong Kong Polytechnic University, Hung Hom, Kowloon, Hong Kong, P.R. China. Email: liqun.qi@polyu.edu.hk. This author's work was supported by the Hong Kong Research Grant Council (Grant No. PolyU 501212, 501913, 15302114 and 15300715).}
}

\date{}

\maketitle

\begin{abstract}
\noindent
In this paper, we introduce a class of variational inequalities, where the involved function is the sum of an arbitrary given vector and a homogeneous polynomial defined by a tensor; and we call it the tensor variational inequality (TVI). The TVI is a natural extension of the affine variational inequality and the tensor complementarity problem. We show that a class of multi-person noncooperative games can be formulated as a TVI. In particular, we investigate the global uniqueness and solvability of the TVI. To this end, we first introduce two classes of structured tensors and discuss some related properties; and then, we show that the TVI has the property of global uniqueness and solvability under some assumptions, which is different from the existed result for the general variational inequality.
\vspace{3mm}

\noindent {\bf Key words:}\hspace{2mm} Tensor variational inequality, global uniqueness and solvability, noncooperative game, strictly positive definite tensor, exceptionally family of elements. \vspace{3mm}

\noindent {\bf AMS subject classifications:}\hspace{2mm} 90C33,
90C30, 65H10. \vspace{3mm}

\end{abstract}

\newpage

\section{Introduction}
\setcounter{equation}{0} \setcounter{Assumption}{0}
\setcounter{Theorem}{0} \setcounter{Proposition}{0}
\setcounter{Corollary}{0} \setcounter{Lemma}{0}
\setcounter{Definition}{0} \setcounter{Remark}{0}
\setcounter{Algorithm}{0}

\hspace{4mm}
The finite dimensional variational inequality (VI) has been studied extensively due to its wide applications in many fields \cite{HP-90,FP-02}. It is called the affine variational inequality if the involved function is linear. The existence and uniqueness of solution to the VI is a basic and important issue in the studies of the VI. It is well known that the VI has at most one solution when the involved function is strictly monotone \cite{HP-90,FP-02,LTY-05}; and a unique solution when the involved function is strongly monotone \cite{HP-90,FP-02}.

It is well known that complementarity problem (CP) is an important subclass of the VIs, which has been studied extensively due to its wide applications \cite{CPS-92,HXQ-06}. Recently, a specific subclass of the CPs, called the tensor complementarity problem (TCP) \cite{SQ-JOTA-15}, has attracted much attention; and many theoretical results about the properties of the solution set of TCP have been developed, including existence of solution \cite{SQ-15,WHB-16,HSW-15r,GLQX-15,YLH-16,SQ-16}, global uniqueness of solution \cite{GLQX-15,BHW-16}, boundedness of solution set \cite{SY-16,SQ-15-2,CQW-15,WHB-16,DLQ-15}, stability of solution \cite{YLH-16}, sparsity of solution \cite{LQX-16}, and so on. In addition, an important application of the TCP was given in \cite{HQ-16}.

Inspired by the development of the TCP, we consider a subclass of the VIs, where the involved function is the sum of an arbitrary given vector and a homogeneous polynomial defined by a tensor; and we call it the tensor variational inequality (TVI). The concerned problem is a natural generalization of the TCP and the affine variational inequality. It is well known that the polynomial optimization problem is an important class of optimization problems, which has been studied extensively \cite{L-01,NDS-06,N-15}. It is easy to see that the TVI is equivalent to a class of polynomial optimization problems. In addition, we show that a class of multi-person noncooperative games can be reformulated as a TVI. These are our motivations to consider the TVI.

In this paper, we mainly investigate the property of global uniqueness and solvability (GUS-property) of the TVI in the case that $0$ belongs to the set involved in the TVI. In this case, we show that there is no strongly monotonously homogeneous polynomial whose degree is larger than 2. In order to investigate the GUS-property of the TVI, we first introduce two classes of structured tensors and discuss some related properties; and then, we show that the TVI has the GUS-property when the involved function is strictly monotone and the involved set contains $0$, which is different from the existed result obtained in the case of the general variational inequality.

The rest of this paper is organized as follows. In Section 2, we recall some basic definitions and results. In Section 3, we introduce the TVI and reformulate a class of multi-person noncooperative games as a TVI. In Section 4, we define two classes of structured tensors and discuss some related properties. In particular, we show that the TVI has the GUS-property under some assumptions. The conclusions are given in Section 5.

\section{Preliminaries}
\setcounter{equation}{0} \setcounter{Assumption}{0}
\setcounter{Theorem}{0} \setcounter{Proposition}{0}
\setcounter{Corollary}{0} \setcounter{Lemma}{0}
\setcounter{Definition}{0} \setcounter{Remark}{0}
\setcounter{Algorithm}{0}

\hspace{4mm}
In this section, we recall some basic concepts and results, which are useful for our subsequent analysis.

Given a nonempty set $X\subseteq \mathbb{R}^n$ and a function $F: X\subseteq \mathbb{R}^n\rightarrow \mathbb{R}^n$, then the variational inequality, denoted by the VI$(X,F)$, is to find a point $x^*\in X$ such that
\begin{equation}\label{vi}
\langle y-x^*, F(x^*)\rangle\geq 0 \mathrm{\quad for\;\; all\;\;} y\in X.
\end{equation}
It is called an affine variational inequality when the function $F$ is linear. Moreover, if the set $X$ is the nonnegative orthant $\mathbb{R}_+^n:=\{x\in \mathbb{R}^n: x\geq 0\}$, then (\ref{vi}) reduces to
$$
x\geq 0,\quad F(x)\geq 0,\quad x^\top F(x)=0,
$$
which is called the complementarity problem, denoted by the CP$(F)$.

In the theoretical studies of the nonlinear variational inequality and complementarity problem, some special types of functions play important roles. The following two classes of functions will be used in this paper.
\begin{Definition}\label{def1} A mapping $F: X\subseteq \mathbb{R}^n\rightarrow \mathbb{R}^n$ is said to be
\begin{itemize}
  \item [(i)]{\bf strictly monotone on $X$} if and only if
  $$\langle F(x)-F(y), x-y\rangle> 0 \mathrm{\quad for\;\; all\;\;} x, y\in X \; \;\mathrm{with}\;\; x\neq y;$$
  \item [(ii)]{\bf strongly monotone on $X$} if and only if there exists a constant $c>0$ such that
  \begin{equation}\label{strongmonotone}\langle F(x)-F(y), x-y\rangle\geq c\|x-y\|^2 \mathrm{\quad for\;\; all\;\;} x, y\in
  X.\end{equation}
\end{itemize}
\end{Definition}

Obviously, a strongly monotone function on $X\subseteq \mathbb{R}^n$ must be strictly monotone on $X$. Moreover, for $X=\mathbb{R}^n$ and an
affine mapping, i.e., $F(x)=Ax+q$, where $A\in \mathbb{R}^{n\times n}$ and $q\in \mathbb{R}^n$, $F$ is strongly monotone if and only if
it is strictly monotone, and if and only if $A$ is positive definite \cite{FP-02}. However, such results do not hold for the general nonlinear function.

The exceptionally family of elements is a powerful tool to investigate the solvability of the VI$(X,F)$ \cite{HHF-04,IZ-00,ZH-99,ZHQ-99,H-03}. There are several different definitions for the exceptionally family of elements. In this paper, we use the following definition.
\begin{Definition}\label{def-efe} \cite[Definition 3.1]{H-03}
Let $\hat{x}\in \mathbb{R}^n$ be an arbitrary given point. A sequence $\{x^r\}_{r>0}$ is said to be an exceptionally family of elements for the VI$(X,F)$ with respect to $\hat{x}$ if the following conditions are satisfied:
\begin{itemize}
  \item $\|x^r\|\rightarrow \infty$ as $r\rightarrow \infty$;
  \item $x^r-\hat{x}\in X$;
  \item there exists $\alpha_r\in (0,1)$ such that, for any $r\geq \|P_X(0)-\hat{x}\|$,
  $$-[F(x^r-\hat{x})+(1-\alpha_r)(x^r-\hat{x})]\in \mathcal{N}_X(x^r-\hat{x}),$$
  where $\mathcal{N}_X(x^r-\hat{x})$ denotes the normal cone of $X$ at $x^r-\hat{x}$ and $P_X(\cdot)$ is the projection operator on $X$.
\end{itemize}
\end{Definition}

The normal cone of $X$ at $x$ is defined by
\begin{equation}\label{normalcone}
  \mathcal{N}_X(x)=\left\{\begin{array}{ll}
  \{z\in \mathbb{R}^n: z^{\top}(y-x)\leq 0, \forall y\in X\},& \mathrm{if}\;\; x\in X,\\
  \varnothing,& \mathrm{otherwise}.
  \end{array}\right.
  \end{equation}
About the relationship between the exceptionally family of elements and the solution of the VI$(X,F)$, we will use the following lemma whose proof can be found in \cite{H-03}.
\begin{Lemma}\cite[Theorem 3.1]{H-03}\label{zeyilemma}
Let $X$ be a nonempty closed convex set in $\mathbb{R}^n$ and $F: X\subseteq \mathbb{R}^n\rightarrow \mathbb{R}^n$ be a continuous function. Then, either the VI$(X,F)$ has a solution or, for any point $\hat{x}\in \mathbb{R}^n$, there exists an exceptionally family of elements for the VI$(X,F)$ with respect to $\hat{x}$.
\end{Lemma}

Throughout this paper, for any positive integer $n$, we use $[n]$ to denote the set $\{1,2,\ldots,n\}$. For any given positive integers $m,r_1,\ldots, r_{m-1}$ and $r_m$, an $m$-order $r_1\times r_2\times\cdots \times r_m$-dimensional real tensor can be denoted by ${\cal A}=(a_{i_1i_2\cdots i_m})$ with $a_{i_1i_2\cdots i_m}\in \mathbb{R}$ for any $i_j\in [r_j]$ and $j\in [m]$. Furthermore, if $r_j=n$ for all $j\in [m]$, then $\cal A$ is called an $m$-order $n$-dimensional real tensor; and we denote the set of all $m$-order $n$-dimensional real tensors by $\mathbb{T}_{m,n}$. In particular, $\mathcal{A}\in \mathbb{T}_{m,n}$ is called a symmetric tensor if the entries $a_{i_1i_2\cdots i_m}$ are invariant under any permutation of their indices. For any ${\cal A}\in \mathbb{T}_{m,n}$ and $x\in \mathbb{R}^n$, $\mathcal{A}x^{m-1}\in \mathbb{R}^n$ is a vector defined by
$$
(\mathcal{A}x^{m-1})_i:=\sum_{i_2,i_3,\cdots,i_m=1}^na_{ii_2\cdots i_m}x_{i_2}x_{i_3}\cdots x_{i_m},\quad \forall i\in [n].
$$

\section{The TVI and an Application}
\setcounter{equation}{0} \setcounter{Assumption}{0}
\setcounter{Theorem}{0} \setcounter{Proposition}{0}
\setcounter{Corollary}{0} \setcounter{Lemma}{0}
\setcounter{Definition}{0} \setcounter{Remark}{0}
\setcounter{Algorithm}{0}

\hspace{4mm}
In this section, we first introduce the TVI and discuss the relationship between it and a class of polynomial optimization problems; and then, give an application of the TVI.

For any ${\cal A}\in \mathbb{T}_{m,n}$, $q\in \mathbb{R}^n$ and a nonempty set $X\subseteq \mathbb{R}^n$, the TVI we considered is given specifically in the following way: Find a vector $x^*\in X$ such that
\begin{equation}\label{tvi}
\langle y-x^*, \mathcal{A}(x^*)^{m-1}+q\rangle\geq 0 \mathrm{\quad for\;\; all\;\;} y\in X,
\end{equation}
which is denoted by the TVI$(X,\mathcal{A},q)$. It should be noted that Song and Qi \cite{SQ-JOTA-15} proposed a TVI$(X,\mathcal{A},q)$ with $q=0$ in a question related to applications of structured tensors; but to the best of our knowledge, the TVI$(X,\mathcal{A},q)$ has not been studied so far even in the case of $q=0$.

The TVI$(X,\mathcal{A},q)$ arises in a natural way in the framework of polynomial optimization problems, which is given as follows:
\begin{Proposition}
For any given symmetric tensor $\mathcal{A}\in \mathbb{T}_{m,n}$ and $q\in \mathbb{R}^n$, let $f(x)=\frac{1}{m}{\cal A}x^m+q^\top x$ be a convex function and $X\subseteq \mathbb{R}^n$ be a nonempty closed convex set. Then, $x^*$ is an optimal solution of the optimization problem
$$
\min\{f(x): x\in X\}
$$
if and only if $x^*$ solves the TVI$(X,\mathcal{A},q)$.
\end{Proposition}
{\bf Proof.} Since $\mathcal{A}$ is symmetric, it follows that $\nabla f(x)=\mathcal{A}x^{m-1}+q$. Then, the result is straightforward from \cite[Page 10]{MS-72}.\ep

In the following, we give an application of the TVI$(X,\mathcal{A},q)$ related to a class of multi-person noncooperative games.

We consider an $m$-person noncooperative game in which each player tries to minimize his own cost. For any $k\in [m]$, let $x^k\in \mathbb{R}^{r_k}$ and $X_k\subseteq \mathbb{R}^{r_k}$ be player $k$'s strategy and strategy set, respectively. We denote
\begin{eqnarray*}
\begin{array}{l}
[m]_{-k}:=[m]\setminus \{k\},\quad n:=\sum\limits_{j\in [m]} r_j,\quad
n_{-k}:=\sum\limits_{j\in [m]_{-k}} r_j,\vspace{2mm}\\
x:=(x^j)_{j\in [m]}\in \mathbb{R}^{r_1}\times\cdots \times\mathbb{R}^{r_m}=\mathbb{R}^n,\vspace{2mm}\\
x^{-k}:=(x^j)_{j\in [m]_{-k}}\in \mathbb{R}^{r_1}\times\cdots \times\mathbb{R}^{r_{k-1}}\times\mathbb{R}^{r_{k+1}}\times\cdots \times\mathbb{R}^{r_{m}}=\mathbb{R}^{n_{-k}},\vspace{2mm}\\
X:=\prod\limits_{j\in [m]} X_j\subseteq \mathbb{R}^{r_1}\times\cdots \times\mathbb{R}^{r_m}=\mathbb{R}^n.
\end{array}
\end{eqnarray*}
For any $k\in [m]$, let $f_k: \mathbb{R}^{r_1}\times\cdots\times \mathbb{R}^{r_m}\rightarrow \mathbb{R}$ denote player $k$'s cost function, which is given by
\begin{eqnarray}\label{func-fk}
f_k(x^k,x^{-k})=\sum_{i_1=1}^{r_1}\sum_{i_2=1}^{r_2}\cdots \sum_{i_m=1}^{r_m}a^k_{i_1i_2\cdots i_m} {x}^{1}_{i_1}{x}^{2}_{i_2}\cdots {x}^{k-1}_{i_{k-1}}x^k_{i_k}{x}^{k+1}_{i_{k+1}}\cdots {x}^{m}_{i_m}.
\end{eqnarray}
Moreover, we use $\mathcal{A}^k=(a^k_{i_1i_2\cdots i_m})$ to denote player $k$'s payoff tensor for any $k\in [m]$.

When the complete information is assumed, for any $k\in [m]$, the $k$th player decides his own strategy by solving the following optimization problem with the opponents' strategy ${x}^{-k}$ fixed:
\begin{eqnarray}\label{pro1}
\min\limits_{x^k}& &f_k(x^k,{x}^{-k})\\\nonumber
\mbox{s.t.}&&x^k\in X_k.
\end{eqnarray}
A tuple ${x^*}:=((x^1)^*,(x^2)^*,\ldots,(x^m)^*)$ satisfying
$$
(x^k)^*\in \mbox{\rm arg}\min_{x^k\in X_k} f_k(x^k,x^{-k}),\quad \forall k\in [m]
$$
is called a Nash equilibrium point of the $m$-person noncooperation game.

In the following, we consider the relationship between the multi-person noncooperation game and the TVI$(X,\mathcal{A},q)$.

\begin{Proposition}\label{gameequtvi}
Suppose that every $X_i\subseteq \mathbb{R}^{r_i}$ is closed and convex, then a tuple
${x^*}:=((x^1)^*,(x^2)^*,\ldots,(x^m)^*)$ is a Nash equilibrium point of the $m$-person noncooperation game if and only if ${x^*}$ is a solution of the VI$(X, F)$ with
\begin{equation}\label{equ30}
F(x)\equiv\left(\nabla_{x^k}f_k(x^k,x^{-k})\right)_{k\in [m]},
\end{equation}
where $\nabla_{x^k}f_k(x^k,x^{-k})$ is the gradient of the function $f_k(x^k,x^{-k})$ defined by (\ref{func-fk}) with respect to $x^k$.
\end{Proposition}
{\bf Proof.} Suppose that $x^*=((x^1)^*,(x^2)^*,\ldots,(x^m)^*)$ is a Nash equilibrium point of the $m$-person noncooperation game, then $(x^k)^*$ is an optimal solution of (\ref{pro1}) for any $k\in [m]$. Since the objective function of the optimization problem (\ref{pro1}) is convex in $x^k$, it follows from the assumption that for any $k\in [m]$, $(x^k)^*$ is an optimal solution of (\ref{pro1}) if and only if
\begin{equation}\label{inequ3}
\left\langle y^k-(x^k)^*, \nabla_{x^k}\left(f_k((x^*)^k,(x^*)^{-k})\right)\right\rangle\geq 0,\quad \forall y^k\in X_k.
\end{equation}
So, ${x^*}$ solves the VI$(X, F)$ with $F$ being defined by (\ref{equ30}).

Conversely, we assume that ${x^*}=((x^1)^*,(x^2)^*,\ldots,(x^m)^*)$ is a solution of the VI$(X, F)$ with $F$ being defined by (\ref{equ30}), then
\begin{equation}\label{equ3}
[F(x^*)]^{\top}(y-x^*)\geq 0,\quad \forall y\in X.
\end{equation}
Due to the arbitrariness of $y$, we let
$$y:=((x^1)^*, (x^2)^*, \ldots, (x^{k-1})^*, y^k, (x^{k+1})^*, \ldots, (x^m)^*)$$
for any $y^k\in X_k$, then (\ref{inequ3}) holds from (\ref{equ3}), which further implies that $(x^k)^*$ is an optimal solution of (\ref{pro1}), i.e., $x^*$ is a Nash equilibrium point of the $m$-person noncooperation game.\ep

In fact, the function $F(x)\equiv(\nabla_{x^k}f_k(x^k,x^{-k}))_{k\in [m]}$ defined in Proposition \ref{gameequtvi} is a homogeneous polynomial function with the degree $m-1$, which can be defined by a tensor. To this end, we first introduce the following symbols: for any tensor $\mathcal{B}\in \mathbb{T}_{m,n}$ and $u^k\in \mathbb{R}^{r_k}$ with $k\in [m]_{-1}$, we denote
$$
\mathcal{B}u^{2}\cdots u^{m}=\left(\begin{array}{c}
\sum\limits_{i_2=1}^{r_2}\cdots\sum\limits_{i_m=1}^{r_m} b_{1i_2\cdots i_n} u^2_{i_2}\cdots u^m_{i_m}\\
\sum\limits_{i_2=1}^{r_2}\cdots\sum\limits_{i_m=1}^{r_m} b_{2i_2\cdots i_n} u^2_{i_2}\cdots u^m_{i_m}\\ \vdots \\ \sum\limits_{i_2=1}^{r_2}\cdots\sum\limits_{i_m=1}^{r_m} b_{r_1i_2\cdots i_n} u^2_{i_2}\cdots u^m_{i_m}
\end{array}\right);
$$
and, for any $k\in [m]$, by using the payoff tensor $\mathcal{A}^k=(a^k_{i_1i_2\cdots i_m})$, we define a new tensor $\bar{\mathcal{A}^k}=(\bar{a}^k_{i_1i_2\cdots i_m})$ with
$$
\bar{a}^k_{i_1i_2\cdots i_m}=a^k_{i_ki_1\cdots i_{k-1}i_{k+1}\cdots i_m}\quad \mbox{\rm for any}\; i_j\in [r_j]\; \mbox{\rm and}\; j\in [m].
$$
Furthermore, we construct a new tensor
$$
\mathcal{A}=(a_{i_1i_2\cdots i_m})\in \mathbb{T}_{m,n},
$$
where for any $i_j\in [n]$ with $j\in [m]$,
\begin{eqnarray*}
a_{i_1i_2\cdots i_m}=\left\{\begin{array}{l}
a^1_{i_1(i_2-r_1)\cdots (i_m-\sum_{j=1}^{m-1}r_j)},\\ \qquad \mbox{\rm if}\;\; i_1\in [r_1], i_2\in [r_1+r_2]\setminus [r_1], \ldots, i_m\in [\sum_{j=1}^mr_j]\setminus [\sum_{j=1}^{m-1}r_j], \vspace{2mm}\\
a^2_{(i_1-r_1)i_2(i_3-r_1-r_2)\cdots (i_m-\sum_{j=1}^{m-1}r_j)},\\ \qquad \mbox{\rm if}\;\; i_1\in [r_1+r_2]\setminus [r_1], i_2\in [r_1], \\
\qquad\quad i_3\in [\sum_{j=1}^3r_j]\setminus [r_1+r_2],\ldots, i_m\in [\sum_{j=1}^mr_j]\setminus [\sum_{j=1}^{m-1}r_j], \vspace{2mm}\\
a^k_{(i_1-\sum_{j=1}^{k-1}r_j)i_2(i_3-r_1)\cdots (i_{k-1}-\sum_{j=1}^{k-3}r_j)i_k(i_{k+1}-\sum_{j+1}^kr_j)\cdots (i_m-\sum_{j=1}^{m-1}r_j)}, \\ \qquad \mbox{\rm if}\;\; k\in [m]\setminus\{1,2\},\; \mbox{\rm and for any given}\; k, i_1\in [\sum_{j=1}^kr_j]\setminus [\sum_{j=1}^{k-1}r_j], \\ \qquad\quad
i_2\in [r_1], i_3\in [r_1+r_2]\setminus [r_1], \ldots, i_{k}\in [\sum_{j=1}^{k-1}r_j]\setminus [\sum_{j=1}^{k-2}r_j],\\ \qquad\quad i_{k+1}\in [\sum_{j=1}^{k+1}r_j]\setminus [\sum_{j=1}^{k}r_j], \ldots, i_m\in [\sum_{j=1}^mr_j]\setminus [\sum_{j=1}^{m-1}r_j],\vspace{2mm}\\
0, \quad \mbox{\rm otherwise}.
\end{array}\right.
\end{eqnarray*}
Then, it is not difficult to get that
\begin{eqnarray}\label{reform}
\mathcal{A}x^{m-1}&=&\left(
\begin{array}{c}
\bar{\mathcal{A}}^1 {x}^{2}\cdots {x}^{m}\\
\vdots\\
\bar{\mathcal{A}}^k {x}^{1}\cdots {x}^{k-1}{x}^{k+1}\cdots {x}^{m}\\
\vdots\\
\bar{\mathcal{A}}^m {x}^{1}{x}^{2}\cdots {x}^{m-1}
\end{array}
\right)=\left(
\begin{array}{c}
\nabla_{x^1}f_1(x^1,x^{-1})\\
\vdots\\
\nabla_{x^k}f_k(x^k,x^{-k})\\
\vdots\\
\nabla_{x^m}f_m(x^m,x^{-m})
\end{array}\right)=F(x).\quad
\end{eqnarray}
Therefore, Proposition \ref{gameequtvi}, together with (\ref{reform}), shows that the concerned $m$-person noncooperative game is to find a Nash equilibrium point $x^*$ satisfying
$$\langle y-x^*, \mathcal{A}(x^*)^{m-1}\rangle\geq 0, \quad \forall y\in X,$$
which is just the TVI$(X,\mathcal{A},q)$ defined by (\ref{tvi}) with $q=0$.

\section{GUS-property of the TVI}
\setcounter{equation}{0} \setcounter{Assumption}{0}
\setcounter{Theorem}{0} \setcounter{Proposition}{0}
\setcounter{Corollary}{0} \setcounter{Lemma}{0}
\setcounter{Definition}{0} \setcounter{Remark}{0}
\setcounter{Algorithm}{0}

\hspace{4mm} The tensor variational inequality (\ref{tvi}) is said to have the GUS-property if it has a unique solution for every $q\in \mathbb{R}^n$. Such an important property has been investigated for variational inequalities \cite{HP-90,FP-02} and complementarity problems \cite{STW-58,MK-77,GS-07,MH-11}. In this section, we discuss the GUS-property of the TVI$(X,\mathcal{A},q)$.

For the general VI, the following results come from \cite{FP-02, HP-90}.
\begin{Lemma}\label{lemma1}Let $X\subseteq \mathbb{R}^n$ be nonempty closed convex and $F: X\rightarrow \mathbb{R}^n$ be continuous.
\begin{itemize}
  \item [(i)]If $F$ is strictly monotone on $X$, then VI$(X,F)$ has at most one solution;
   \item [(ii)]If $F$ is strongly monotone on $X$, then VI$(X,F)$ has a unique solution.
   \end{itemize}
\end{Lemma}

Let $F: X\subseteq \mathbb{R}^n\rightarrow \mathbb{R}^n$ be defined by
\begin{eqnarray}\label{func-f-1}
F(x):=\mathcal{A}x^{m-1}+q,
\end{eqnarray}
where $\mathcal{A}\in \mathbb{T}_{m,n}$ with $m>2$ and $q\in \mathbb{R}^n$. Then, we have the following observation.
\begin{Proposition}\label{prop-1}
For any tensor $\mathcal{A}\in \mathbb{T}_{m,n}$ with $m>2$ and $q\in \mathbb{R}^n$, let the function $F$ be defined by (\ref{func-f-1}).
Suppose that $0\in X\subseteq \mathbb{R}^n$, then the function $F$ is not strongly monotone on $X$.
\end{Proposition}
{\bf Proof.} Suppose that there exist a vector $q\in \mathbb{R}^n$ and a tensor $\mathcal{A}\in \mathbb{T}_{m,n}$ with $m>2$ such that the function $F$ defined by (\ref{func-f-1}) is strongly monotone on $X$, then there exists a positive constant $c$ such that (\ref{strongmonotone}) holds for any $x, y\in X$. Let $y=0\in X$, then we get from (\ref{strongmonotone}) that
\begin{equation}\label{inequ4}
\mathcal{A}x^{m}\geq c \|x\|^2 \quad \textrm{for any}\; x\in X.
\end{equation}
For any $x\neq 0$, it follows from (\ref{inequ4}) that
\begin{equation}\label{inequ5}
\mathcal{A}\left(\frac{x}{\|x\|}\right)^{m}\geq c \left\|\left(\frac{x}{\|x\|}\right)\right\|^2\|x\|^{2-m}.
\end{equation}
Since $\left\|\frac{x}{\|x\|}\right\|=1$, it follows that the left-hand side of the inequality (\ref{inequ5}) is bounded;  but when $\|x\|\rightarrow 0$, it is obvious that the right-hand side of the inequality (\ref{inequ5}) tends to $\infty$, which leads to a contradiction. Therefore, there exists no the strongly monotone function $F$ in the form of $\mathcal{A}x^{m-1}+q$ for any $q\in \mathbb{R}^n$ and $\mathcal{A}\in \mathbb{T}_{m,n}$ with $m>2$.\ep

From Lemma \ref{lemma1} (ii) and Proposition \ref{prop-1}, a natural question is \emph{whether or not the VI$(X,F)$ has the GUS-property when $0\in X$ and the function $F$ is defined by (\ref{func-f-1}) where $\mathcal{A}\in \mathbb{T}_{m,n}$ with $m>2$ and $q\in \mathbb{R}^n$}. In this section, we answer this question. To this end, we firstly introduce two new classes of tensors in the next
subsection and discuss the relationship between them.

\subsection{Relationship of Two Classes of Tensors}
\hspace{4mm} In this subsection, we introduce two new classes of structured tensors and discuss the relationship between them.
\begin{Definition}\label{def3}
Given a nonempty set $X\subseteq \mathbb{R}^n$. A tensor $\mathcal{\mathcal{A}}\in \mathbb{T}_{m,n}$ is said to be
\begin{itemize}
  \item [(i)]{\bf positive definite on $X$}  if and only if $\mathcal{A}x^{m}> 0$ for any $x\in X$ and $x\neq 0$, and
  \item [(ii)]{\bf strictly positive definite on $X$}  if and only if
  $$(x-y)^\top(\mathcal{A}x^{m-1}-\mathcal{A}y^{m-1})> 0\quad \mbox{\rm for any}\; x,y\in X\;\mbox{\rm with}\; x\neq y.$$
\end{itemize}
$\mathcal{\mathcal{A}}\in \mathbb{T}_{m,n}$ is said to be a strictly positive definite tensor if it is strictly positive definite on $\mathbb{R}^n$.
\end{Definition}

When $X=\mathbb{R}^n$, the positive definite tensor on $X$ defined by Definition \ref{def3} (i) is just the positive definite tensor
defined in \cite{Q-05}; and when $X=\mathbb{R}^n_+$, the positive definite tensor on $X$ defined by Definition \ref{def3} (i) is just the strictly copositive tensor defined in \cite{Q-13}. From Definitions \ref{def1} and \ref{def3}, it is easy to see that the function $F$ defined by (\ref{func-f-1}) is strictly monotone on $X$ if and only if the tensor $\cal A$ is strictly positive definite on $X$.

A basic question is \emph{whether or not there exists a strictly positive definite tensor on some subset of $\mathbb{R}^n$}. The following example gives a positive answer to this question.

\begin{Example}\label{exam1}
Let $\mathcal{A}=(a_{ijkl})\in \mathbb{T}_{4,2}$, where $a_{1111}=a_{2222}=1$,  and the others equal to zero. Then, $\mathcal{A}$ is a strictly positive definite tensor on any subset $X$ of $\mathbb{R}^2$.
\end{Example}

It only needs to prove that $\mathcal{A}$ is strictly positive definite on $\mathbb{R}^2$.

Since
$$
\mathcal{A}x^3=\left(\begin{array}{c}x_1^3\\x_2^3\end{array}\right),
$$
it follows that for any $x, y\in \mathbb{R}^2$,
\begin{eqnarray}\label{equ4}
(x_1-y_1)[(\mathcal{A}x^3)_1-(\mathcal{A}y^3)_1]&=&(x_1-y_1)(x_1^3-y_1^3)\nonumber\\
&=&(x_1-y_1)^2(x_1^2+x_1y_1+y_1^2);
\end{eqnarray}
\begin{eqnarray}\label{equ5}
(x_2-y_2)[(\mathcal{A}x^3)_2-(\mathcal{A}y^3)_2]&=&(x_2-y_2)(x_2^3-y_2^3)\nonumber\\
&=&(x_2-y_2)^2(x_2^2+x_2y_2+y_2^2).
\end{eqnarray}

For any $s,t\in \mathbb{R}$, we discuss the following three cases.
\begin{itemize}
\item [(I)] $|s|\neq |t|$. In this case, we have
  $$s^2+st+t^2> 2|s||t|+st=\left\{\begin{array}{rll}3st\geq 0,& \mbox{if} &st\geq 0,\\
  -st>0,& \mbox{if} &st<0, \end{array}\right.$$
  which implies that $s^2+st+t^2>0$.
\item [(II)] $s=t$. In this case, we have
  $$(s-t)^2(s^2+st+t^2)=0.$$
\item [(III)] $s=-t\neq 0$. In this case, we have
  $$(s-t)^2(s^2+st+t^2)=4s^4>0.$$
\end{itemize}

Now, for any $x,y\in \mathbb{R}^2$ and $x\neq y$, it follows that either $x_1\neq y_1$ or $x_2\neq y_2$. Therefore, by combining cases
(I)-(III) with (\ref{equ4}) and (\ref{equ5}) we have
$$
(x-y)^\top(\mathcal{A}x^3-\mathcal{A}y^3)=\sum_{i=1}^2(x_i-y_i)^2(x_i^2+x_iy_i+y_i^2)>0,
$$
which demonstrates that $\mathcal{A}$ is a strictly positive
definite tensor on $\mathbb{R}^2$.\ep

In the following, we discuss the relationship between two classes of tensors defined by Definition \ref{def3}.

\begin{Proposition}\label{proconclusions}
Suppose that $0\in X\subseteq \mathbb{R}^n$. Then, a strictly positive definite tensor on $X$ must be positive definite on $X$.
\end{Proposition}
{\bf Proof.} Given a tensor $\mathcal{A}\in \mathbb{T}_{m,n}$. Take $y=0\in X$, it follows from Definition \ref{def3}(ii) that for any $x\in X$ with $x\neq 0$,
$$
\mathcal{A}x^m=(x-0)^\top\left(\mathcal{A}x^{m-1}-\mathcal{A}0^{m-1}\right)>0,
$$
which, together with Definition \ref{def3}(i), implies that $\mathcal{A}$ is positive definite on $X$.  \ep

However, if $m>2$, a positive definite tensor on $X$ is not necessary a strictly positive definite tensor on $X$, which can be seen in the following example.

\begin{Example}
Denote $X:=\mathbb{R}^2_+$. Let $\mathcal{A}=(a_{ijkl})\in \mathbb{T}_{4,2}$, where $a_{1111}=a_{2222}=a_{2112}=1$, $a_{1122}=-1$, and the others equal to zero. Then, $\mathcal{A}$ is positive definite on $X$ but not
strictly positive definite on $X$.
\end{Example}

Firstly, we show that $\mathcal{A}$ is positive definite on $X$. Since
$$
\mathcal{A}x^3=\left(\begin{array}{c}x_1^3-x_1x_2^2\\x_2^3+x_1^2x_2\end{array}\right),
$$
it follows that for any $x\in \mathbb{R}^2\setminus \{0\}$,
$$
x^\top\mathcal{A}x^3=x_1^4-x_1^2x_2^2+x_2^4+x_1^2x_2^2=x_1^4+x_2^4>0.
$$
Hence, $\mathcal{A}$ is positive definite on $\mathbb{R}^2$. Of course, $\mathcal{A}$ is positive definite on $X$.

Secondly, we show that $\mathcal{A}$ is not a strictly positive definite tensor on $X$. To this end, for any $\mu\in \mathbb{R}_+$ with
$\mu\neq 0$, let $x=(2\mu, 3\mu)^\top$ and $y=(\mu, 3\mu)^\top$, then
$x,y\in X$, $x\neq y$ and
\begin{eqnarray*}
(x-y)^\top(\mathcal{A}x^3-\mathcal{A}y^3)
&=&(x_1-y_1)[(\mathcal{A}x^3)_1-(\mathcal{A}y^3)_1]+(x_2-y_2)[(\mathcal{A}x^3)_2-(\mathcal{A}y^3)_2]\\
&=&(2\mu-\mu)[(2\mu)^3-2\mu(3\mu)^2-(\mu^3-\mu(3\mu)^2)]+0\\
&=&-2\mu^4\\
&<&0.
\end{eqnarray*}
Therefore, $\mathcal{A}$ is not strictly positive definite on $X$. \ep

\subsection{Uniqueness of Solution to the TVI}

\hspace{4mm} In this subsection, we investigate the GUS-property of the TVI$(X,\mathcal{A},q)$.
\begin{Theorem}\label{atmostonesolution}
Let $X\subseteq \mathbb{R}^n$ be a nonempty closed convex set and $\mathcal{A}\in \mathbb{T}_{m,n}$ be a strictly positive definite tensor on $X$. Then, for any given $q\in \mathbb{R}^n$, the TVI$(X,\mathcal{A},q)$ has at most one solution.
\end{Theorem}
{\bf Proof.} Since $\mathcal{A}$ is a strictly positive definite tensor on $X$, it follows from Definition \ref{def3} (ii) that the function $\mathcal{A}x^{m-1}+q$ is strictly monotone on $X$ for any $q\in \mathbb{R}^n$. So, the desired result holds
from Lemma \ref{lemma1} (i). \ep

\begin{Theorem}\label{nonemptysolution}
Let $X\subseteq \mathbb{R}^n$ be a nonempty closed convex set with $0\in X$ and $\mathcal{A}\in \mathbb{T}_{m,n}$ be a positive definite tensor on $X$. Then, for any given $q\in \mathbb{R}^n$, the solution set of the TVI$(X,\mathcal{A},q)$ is nonempty and compact.
\end{Theorem}
{\bf Proof.} If the set $X$ is bounded, then the result is obvious from \cite{HP-90,HS-66}. In the following, we assume that the set $X$ is unbounded.

Suppose that the TVI$(X,\mathcal{A},q)$ has no solution, then for $\hat{x}=0\in \mathbb{R}^n$, it follows from Lemma \ref{zeyilemma} that there exists an exceptionally family of elements $\{x^r\}_{r>0}$ for the TVI$(X,\mathcal{A},q)$ with respect
to $0$. That is, we have
\begin{itemize}
  \item [(a)]$\|x^r\|\rightarrow \infty$ as $r\rightarrow \infty$;
  \item [(b)]$x^r\in X$ for any positive integer $r$;
  \item [(c)]there exists $\alpha_r\in (0,1)$ such that, for any $r\geq \|P_X(0)\|$,
  $$
  -[\mathcal{A}(x^r)^{m-1}+(1-\alpha_r)x^r]\in \mathcal{N}_X(x^r).
  $$
\end{itemize}
From the above (c) and the definition of the normal cone, we have
  \begin{eqnarray*}
  [\mathcal{A}(x^r)^{m-1}+(1-\alpha_r)x^r]^\top(y-x^r)\geq 0\quad \mbox{\rm for any}\; y\in X,
  \end{eqnarray*}
which can be rewritten as
  \begin{equation}\label{normal1}
  [\mathcal{A}(x^r)^{m-1}]^\top(y-x^r)\geq (\alpha_r-1)(x^r)^\top(y-x^r)\quad \mbox{\rm for any}\; y\in X,
  \end{equation}
From the above (a), it holds that $\|x^r\|>0$ for sufficiently large $r$. So, by dividing $\|x^r\|^{m}$ in both sides of the inequality (\ref{normal1}), we get
\begin{eqnarray*}
\left[\mathcal{A}\frac{(x^r)^{m-1}}{\|x^r\|^{m-1}}\right]^\top\left(\frac{y}{\|x^r\|}-\frac{x^r}{\|x^r\|}\right)\geq
\frac{\alpha_r-1}{\|x^r\|^{m-2}}\left(\frac{x^r}{\|x^r\|}\right)^\top\left(\frac{y}{\|x^r\|}-\frac{x^r}{\|x^r\|}\right).
\end{eqnarray*}
Let $z^r=\frac{x^r}{\|x^r\|}$, then the above inequality becomes
\begin{equation}\label{theorm-5.2-1}
[\mathcal{A}(z^r)^{m-1}]^\top\left(\frac{y}{\|x^r\|}-z^r\right)\geq
\frac{\alpha_r-1}{\|x^r\|^{m-2}}(z^r)^\top\left(\frac{y}{\|x^r\|}-z^r\right).
\end{equation}
Since the sequence $\{z^r\}$ is bounded, there exists a convergent subsequence. Without lose of generality, we denote this subsequence
by $\{z^r\}$ and its limit point by $z^*$. Noting that $\alpha_r\in(0,1)$ and $y\in X$ is an arbitrary given vector, by letting $r\rightarrow \infty$, it follows from (\ref{theorm-5.2-1}) that $[\mathcal{A}(z^*)^{m-1}]^\top(-z^*)\geq 0$, i.e.,
  \begin{eqnarray}\label{theorm-5.2-2}
  \mathcal{A}(z^*)^{m}\leq 0.
  \end{eqnarray}
Next, we show that $z^*\in X$. Since $\|x^r\|\rightarrow \infty$ as $r\rightarrow \infty$, it follows that $\frac{1}{\|x^r\|}<1$ with sufficiently large $r$. Furthermore, since $0\in X$ and $X$ is convex, it follows from the above (b) that for sufficiently large $r$,
$$
z^r=\frac{x^r}{\|x^r\|}=\left(1-\frac{1}{\|x^r\|}\right)0+\frac{1}{\|x^r\|}x^r\in X.
$$
Thus, by the fact that the set $X$ is closed, we get
$$
z^*\in X.
$$
This, together with (\ref{theorm-5.2-2}), contradicts that $\mathcal{A}$ is a positive definite tensor on $X$. Therefore, the
TVI$(X,\mathcal{A},q)$ has at least one solution when $\mathcal{A}$ is a positive definite tensor on $X$.

Denote the solution set of the TVI$(X,\mathcal{A},q)$ by SOL$(X,\mathcal{A},q)$. Suppose that the sequence $\{x^k\}\subseteq$SOL$(X,\mathcal{A},q)$ and $x^k\rightarrow x^*$ as $k\rightarrow \infty$, then it follows that
$$
(y-x^k)^\top\left[\mathcal{A}(x^k)^{m-1}+q\right]\geq 0 \quad \mbox{\rm for all}\; y\in X.$$
Thus, let $k\rightarrow \infty$, we get
$$
(y-x^*)^\top\left[\mathcal{A}(x^*)^{m-1}+q\right]\geq 0 \quad \mbox{\rm for all}\; y\in X.$$
That is, $x^*\in$SOL$(X,\mathcal{A},q)$. So, the solution set of the TVI$(X,\mathcal{A},q)$ is closed.

Suppose that the solution set of the TVI$(X,\mathcal{A},q)$ is unbounded, then there exists a sequence $\{x^k\}\subseteq$SOL$(X,\mathcal{A},q)$ such that
$\|x^k\|\rightarrow \infty$ as $k\rightarrow \infty$. Since
$$
(y-x^k)^\top\left[\mathcal{A}(x^k)^{m-1}+q\right]\geq 0 \quad \mbox{\rm for all}\; y\in X,$$
which leads to
$$
\left(\frac{y}{\|x^k\|}-\frac{x^k}{\|x^k\|}\right)^\top\left[\mathcal{A}\left(\frac{x^k}{\|x^k\|}\right)^{m-1}+\frac{q}{\|x^k\|^{m-1}}\right]\geq 0.
$$
Let $k\rightarrow \infty$ and denote $x^*=\lim_{k\rightarrow
\infty}\frac{x^k}{\|x^k\|}$, then we have that
$$
x^*\in X\quad \mbox{\rm and}\quad -\mathcal{A}({x}^*)^{m}\geq 0,
$$
which contradicts the condition that $\mathcal{A}$ is a positive definite tensor on $X$. So, the solution set of the TVI$(X,\mathcal{A},q)$ is bounded.

The proof is complete. \ep

\begin{Corollary}\label{corollary3}
Let $X\subseteq \mathbb{R}^n$ be a nonempty closed convex set with $0\in X$ and $\mathcal{A}\in \mathbb{T}_{m,n}$ be a strictly positive definite tensor on $X$. Then, for any given $q\in \mathbb{R}^n$, the solution set of the TVI$(X,\mathcal{A},q)$ is nonempty and compact.
\end{Corollary}
{\bf Proof.} Since $0\in X$, it follows from Proposition \ref{proconclusions} that a strictly positive definite tensor on $X$ is necessary a positive definite tensor on $X$. Thus, the result is obvious. \ep

\begin{Theorem}\label{theorem-1}
Let $X\subseteq \mathbb{R}^n$ be a nonempty closed convex set with $0\in X$ and $\mathcal{A}\in \mathbb{T}_{m,n}$ be a strictly positive definite tensor on $X$. Then, for any given $q\in \mathbb{R}^n$, the TVI$(X,\mathcal{A},q)$ has a unique solution.
\end{Theorem}
{\bf Proof.} By virtue of Theorem \ref{atmostonesolution} and Corollary \ref{corollary3}, the result is straightforward. \ep

Equivalently, we have the following result.
\begin{Corollary}
Let $X\subseteq \mathbb{R}^n$ be a nonempty closed convex set with $0\in X$ and $\mathcal{A}\in \mathbb{T}_{m,n}$.
Suppose that the function $F(x):=\mathcal{A}x^{m-1}+q$ is strictly monotone on $X$, then the VI$(X,F)$ has a unique solution for any $q\in \mathbb{R}^n$.
\end{Corollary}

Let $X\subseteq \mathbb{R}^n$ be a nonempty closed convex set and the function $F$ be given by $F(x)=\mathcal{A}x^{m-1}+q$ where $\mathcal{A}\in \mathbb{T}_{m,n}$ and $q\in \mathbb{R}^n$. We have showed that, in the case of $0\in X$, the VI$(X,F)$ has the GUS-property if the function $F$ is strictly monotone on $X$. What would happen if $0\notin X$? From Lemma \ref{lemma1}, we know that the VI$(X,F)$ has the GUS-property if the function $F$ is strongly monotone on $X$. A natural question is \emph{whether or not there exists a strongly monotone function $F(x)=\mathcal{A}x^{m-1}+q$ (with $m>2$) on $X$ with $0\notin X$}. The following example gives a positive answer to this question.

\begin{Example}
Let
\begin{equation}\label{closet}
X:=\{(u,1)^\top:u\in \mathbb{R}, u\geq1\},
\end{equation}
and $\mathcal{A}\in \mathbb{T}_{m,n}$ be defined in Example \ref{exam1}, then $F(x):=\mathcal{A}x^{m-1}+q$ with any $q\in \mathbb{R}^2$ is strongly monotone on $X$.
\end{Example}

For any $x,y\in X$, it follows that there exist $u\geq 1$ and $v\geq 1$ such that $x=(u,1)^\top$ and $y=(v,1)^\top$. Furthermore, for any $q\in \mathbb{R}^2$, we have
$$
(x-y)^\top\left[F(x)-F(y)\right]=(x-y)^\top\left(\mathcal{A}x^3-\mathcal{A}y^3\right)=(u-v)^2(u^2+uv+v^2);
$$
but for $\mu=1$, we have
$$
\mu\|x-y\|^2=(u-v)^2.
$$
Obviously,
$$
(u-v)^2(u^2+uv+v^2)\geq 3uv(u-v)^2\geq 3(u-v)^2\geq (u-v)^2.
$$
Thus, for any $x,y\in X$ and $q\in \mathbb{R}^2$, there exists a constant $\mu=1$ such that
$$
(x-y)^\top\left(\mathcal{A}x^3-\mathcal{A}y^3\right)\geq \mu\|x-y\|^2.
$$
So, the function $F$ is strongly monotone on the set $X$ defined by (\ref{closet}). \ep

Therefore, when $X\subseteq \mathbb{R}^n$ is a nonempty closed convex set with $0\notin X$, from Lemma \ref{lemma1} (ii), we
know that the TVI$(X,\mathcal{A},q)$ has a unique solution on $X$ if the function $\mathcal{A}x^{m-1}+q$ is strongly monotone on $X$. We do not know whether the condition of strong monotonicity can be weaken or not in this case.

Before the end of this section, we illustrate that a strictly monotone function $\mathcal{A}x^{m-1}+q$ on $X\subseteq \mathbb{R}^n$ is not necessarily strongly monotone on $X$ when $0\notin X$.
\begin{Example}
Let $\mathcal{A}\in \mathbb{T}_{m,n}$ be defined in Example \ref{exam1} and $X:=\{(u,1)^\top: u\in \mathbb{R}\}$. Then, for any $q\in \mathbb{R}^2$, the function $\mathcal{A}x^{m-1}+q$ is strictly monotone on $X$ but not strongly monotone on $X$.
\end{Example}

Firstly, from Example \ref{exam1}, it is obvious that the tensor $\mathcal{A}$ is strictly positive definite on $X$. Therefore,
the function $\mathcal{A}x^{m-1}+q$ is strictly monotone on $X$.

Secondly, we show that the function $\mathcal{A}x^{m-1}+q$ is not strongly monotone on $X$. Suppose that $\mathcal{A}x^{m-1}+q$ is strongly
monotone on $X$, then there exists a scalar $\mu_0>0$ such that
\begin{equation}\label{strongmo}
(x-y)^\top\left(\mathcal{A}x^3-\mathcal{A}y^3\right)\geq\mu_0\|x-y\|^2\quad \mbox{\rm for any}\; x,y\in X.
\end{equation}

Now, take $x^0=(\sqrt{\mu_0}, 1)^\top\in X$ and $y^0=(-\frac{\sqrt{\mu_0}}{2}, 1)^\top\in X$, then
\begin{eqnarray*}
\left(x^0-y^0\right)^\top\left[\mathcal{A}(x^0)^3-\mathcal{A}(y^0)^3\right]
&=&\left(x_1^0-y_1^0\right)^2\left[(x_1^0)^2+x_1^0y_1^0+(y_1^0)^2\right]\\
&=&\left[\sqrt{\mu_0}+\frac{\sqrt{\mu_0}}{2}\right]^2\left[(\sqrt{\mu_0})^2-\sqrt{\mu_0}\cdot\frac{\sqrt{\mu_0}}{2}+\left(\frac{\sqrt{\mu_0}}{2}\right)^2\right]\\
&=&\left(\frac{3\sqrt{\mu_0}}{2}\right)^2\left(\mu_0-\frac{\mu_0}{2}+\frac{\mu_0}{4}\right)=\frac{27}{16}\mu_0^2
\end{eqnarray*}
and
$$
\mu_0\|x^0-y^0\|^2=\mu_0\left[(x_1^0-y_1^0)^2+(x_2^0-y_2^0)^2\right]
=\mu_0\left[\sqrt{\mu_0}+\frac{\sqrt{\mu_0}}{2}\right]^2=\frac{9}{4}\mu_0^2.
$$
These yield that
$$
(x^0-y^0)^\top\left[\mathcal{A}(x^0)^3-\mathcal{A}(y^0)^3\right]<\mu_0\|x^0-y^0\|^2,
$$
which contradicts the inequality (\ref{strongmo}). So, the function $\mathcal{A}x^{m-1}+q$ is not strongly monotone on $X$.\ep

\section{Conclusions}
\setcounter{equation}{0} \setcounter{Assumption}{0}
\setcounter{Theorem}{0} \setcounter{Proposition}{0}
\setcounter{Corollary}{0} \setcounter{Lemma}{0}
\setcounter{Definition}{0} \setcounter{Remark}{0}
\setcounter{Algorithm}{0}

\hspace{4mm} In this paper, we introduced the tensor variational inequality which arises in a natural way in the framework of
polynomial optimization problems when the involved tensor is symmetric; and showed that a class of multi-person noncooperation games can be reformulated as a class of tensor variational inequalities. In particular, we showed that the tensor variational inequality TVI$(X,\mathcal{A},q)$ has the GUS-property when the function $F(x):=\mathcal{A}x^{m-1}+q$ is strictly monotone on $X$ and $0\in X$, which is different from the existed result obtained in the general variational inequality.

It should be pointed out that we have just done some initial research for the tensor variational inequality in this paper. Many questions need to be answered in the future. Here, we provide three questions as follows.
\begin{Question}
For the TVI$(X,\mathcal{A},q)$ with $X$ being a nonempty closed convex set, when $0\in X$, we showed that the TVI$(X,\mathcal{A},q)$ has the GUS-property if $F(x):=\mathcal{A}x^{m-1}+q$ is strictly monotone on $X$. It is worth investigating whether the condition $0\in X$ can be removed or weaken or not.
\end{Question}
\begin{Question}
How to design effective algorithms to solve the TVI$(X,\mathcal{A},q)$ by using the specific structure of the tensor
$\mathcal{A}$?
\end{Question}
\begin{Question}
In \cite{G-16}, the author investigated the properties of the general polynomial complementarity problem denoted by the PCP$(f)$ with
\begin{eqnarray}\label{e-final}
f(x)=\mathcal{A}_mx^{m-1}+\mathcal{A}_{m-1}x^{m-2}+\cdots+\mathcal{A}_2x+\mathcal{A}_1,
\end{eqnarray}
where $\mathcal{A}_k$ is a tensor of order $k$ and $\mathcal{A}_kx^{k-1}$ is a polynomial mapping for any $k\in [m]$. If we use the polynomial function $f$ defined by (\ref{e-final}) to replace the function $\mathcal{A}x^{m-1}+q$ in the TVI$(X,\mathcal{A},q)$, i.e., find a vector $x^*\in X$ such that
\begin{eqnarray*}
\langle y-x^*, f(x^*)\rangle\geq 0 \mathrm{\quad for\;\; all\;\;} y\in X,
\end{eqnarray*}
then we call it the polynomial variational inequality, denoted by the PVI$(X,f)$.  What are the properties of solution to the PVI$(X,f)$?
\end{Question}



\end{document}